\documentclass[titlepage,draft,12pt]{article} 
\usepackage{amsfonts,latexsym,pstricks,pst-plot} 
\usepackage[a4paper]{geometry}

\newcommand{\ocont}{$\omega$-continuous}
\newcommand{\G}{\mathcal{G}}
\newcommand{\trebar}[1]{\|#1\|}
\newcommand{\lk}{\lambda_{k}}
\newcommand{\ep}{\varepsilon}
\newcommand{\re}{{\mathbb{R}}}
\newcommand{\n}{{\mathbb{N}}}

\newcommand{\vs}{\vspace{.2cm}}
\newcommand{\prf}{{\sc Proof.}\ }
\newcommand{\qed}{{\penalty 10000\mbox{$\quad\Box$}}\bigskip}

\newcommand{\m}[1]{m(\au{#1}^{2})}
\newcommand{\au}[1]{|A^{1/2}#1|}

\newcommand{\cep}{c_{\ep}}

\newtheorem{thm}{Theorem}[section]
\newtheorem{thmbibl}{Theorem}
\newtheorem{rmk}[thm]{Remark}

\newtheorem{lemma}[thm]{Lemma}
\newtheorem{open}[thm]{Open problem}

\title{A uniqueness result for Kirchhoff equations with non-Lipschitz
nonlinear term}
\author{Marina Ghisi\vs\\ {\normalsize
Universit\`a degli Studi di Pisa} \\{\normalsize Dipartimento di
Matematica ``Leonida Tonelli''}\\
{\normalsize 
PISA (Italy)}\\  
{\normalsize e-mail: \texttt{ghisi@dm.unipi.it}}\and 
Massimo Gobbino\vs\\ {\normalsize Universit\`a degli Studi di Pisa} 
\\{\normalsize Dipartimento di Matematica Applicata ``Ulisse Dini''}\\ 
{\normalsize 
 PISA (Italy)}\\  
{\normalsize e-mail: \texttt{m.gobbino@dma.unipi.it}}}
\date{}

\begin{document}
\maketitle
\begin{abstract}
	We consider the second order Cauchy problem 
	$$u''+\m{u}Au=0, 
	\hspace{2em}
	u(0)=u_{0},\quad
	u'(0)=u_{1},$$
	where $m:[0,+\infty)\to[0,+\infty)$ is a continuous function, and
	$A$ is a self-adjoint nonnegative operator with dense domain on a
	Hilbert space.
	
	It is well known that this problem admits local-in-time solutions
	provided that $u_{0}$ and $u_{1}$ are regular enough, depending on
	the continuity modulus of $m$. It is also well known that the
	solution is unique when $m$ is locally Lipschitz continuous.
	
	In this paper we prove that if either $\langle
	Au_{0},u_{1}\rangle\neq 0$, or
	$|A^{1/2}u_{1}|^{2}\neq\m{u_{0}}|Au_{0}|^{2}$, then the local
	solution is unique even if $m$ is not Lipschitz continuous.
		
\vspace{1cm}

\noindent{\bf Mathematics Subject Classification 2000 (MSC2000):}
35L70, 35L80, 35L90.

\vspace{1cm} 

\noindent{\bf Key words:} uniqueness, integro-differential hyperbolic
equation, continuity modulus, Kirchhoff equations, Gevrey spaces.
\end{abstract}
 
\section{Introduction}

Let $H$ be a real Hilbert space.  For every $x$ and $y$ in $H$, let
$|x|$ denote the norm of $x$, and let $\langle x,y\rangle$ denote the
scalar product of $x$ and $y$.  Let $A$ be an unbounded linear
operator on $H$ with dense domain $D(A)$.  We always assume that $A$
is self-adjoint and nonnegative, so that the power $A^{\alpha}$ is
defined for every $\alpha\geq 0$ in a suitable domain $D(A^{\alpha})$.

Given a continuous function $m:[0,+\infty)\to[0,+\infty)$ we consider 
the Cauchy problem
\begin{equation}
	u''(t)+\m{u(t)}Au(t)=0, 
	\hspace{2em}\forall t\in[0,T),
	\label{pbm:h-eq}
\end{equation}
\begin{equation}
	u(0)=u_0,\hspace{3em}u'(0)=u_1.
	\label{pbm:h-data}
\end{equation}
 
It is well known that (\ref{pbm:h-eq}), (\ref{pbm:h-data}) is the
abstract setting of the Cauchy-boundary value problem for the
quasilinear hyperbolic integro-differential partial differential
equation
\begin{equation}
	u_{tt}(t,x)-
	m{\left(\int_{\Omega}\left|\nabla u(t,x)\right|^2\,dx\right)}
	\Delta u(t,x)=0
	\hspace{2em}
	\forall(x,t)\in\Omega\times[0,T),
	\label{eq:k}
\end{equation}
where $\Omega\subseteq\re^{n}$ is an open set, and $\nabla u$ and
$\Delta u$ denote the gradient and the Laplacian of $u$ with respect
to the space variables.

A lot of papers have been devoted to existence of local or global
solutions to (\ref{pbm:h-eq}), (\ref{pbm:h-data}).  The interested
reader is referred to the references quoted in \cite{ap} or in the
more recent papers \cite{gg:k-derloss}, \cite{hirosawa2},
\cite{manfrin1}.

In particular a local-in-time solution to (\ref{pbm:h-eq}),
(\ref{pbm:h-data}) is known to exist provided that the initial data
$u_{0}$ and $u_{1}$ are regular enough.  As in the linear case, the
required regularity depends on the continuity modulus $\omega$ of $m$,
and on the strict hyperbolicity ($m(\sigma)\geq\nu>0$ for every
$\sigma\geq 0$) or weakly hyperbolicity ($m(\sigma)\geq0$ for every
$\sigma\geq 0$) of equation (\ref{pbm:h-eq}).  A rough sketch of the
situation for the strictly hyperbolic case is provided by the
following scheme:
		\begin{eqnarray*}
			\omega(\sigma)=o(\sigma) & \to & \mbox{analytic data,}  \\
			\omega(\sigma)=\sigma^{\alpha}\mbox{ (with $\alpha\in(0,1)$)} & \to
			& \mbox{Gevrey space
			$\G_{s}(A)$ with $s=(1-\alpha)^{-1}$,}  \\
			\omega(\sigma)=\sigma|\log \sigma| & \to & 
			\mbox{$D(A^{\infty})$ (finite derivative loss),}  \\
			\omega(\sigma)=\sigma & \to & \mbox{$D(A^{3/4})\times
			D(A^{1/4})$ (no derivative loss)}.
		\end{eqnarray*}
	
More regularity is required in the weakly hyperbolic case, according
to the following scheme:
		\begin{eqnarray*}
			\omega(\sigma)=o(\sigma) & \to & \mbox{analytic data,}  \\
			\omega(\sigma)=\sigma^{\alpha}\mbox{ (with $\alpha\in(0,1)$)} & \to
			& \mbox{Gevrey space
			$\G_{s}(A)$ with $s=1+\alpha/2$,}  \\
			\omega(\sigma)=\sigma & \to & \mbox{Gevrey space
			$\G_{3/2}(A)$}.
		\end{eqnarray*}

We refer to section~\ref{sec:statements} for a formal statement
(Theorem~\ref{thm:hirosawa}), and for precise definitions of the
functional spaces in the abstract setting.

In this paper we focus on the uniqueness problem for these local
solutions.

It is well known that uniqueness holds whenever $m$ is Lipschitz
continuous.  This result has been proved for example in~\cite{ap} in
the strictly hyperbolic case with initial data in $D(A^{3/4})\times
D(A^{1/4})$, and in~\cite{as} in the weakly hyperbolic case with
analytic initial data.  In the weakly hyperbolic case the same
argument can be easily extended to initial data in the Gevrey class
$\G_{3/2}(A)$, which is the largest space where local existence can be
proved (of course in the weakly hyperbolic case with a Lipschitz
continuous $m$).

When $m$ is not Lipschitz continuous the uniqueness problem seems to
be widely unexplored.  To our knowledge indeed this case has been
considered only in section~4 of~\cite{as}, where two results are
presented.  The first one is a one-dimensional example where problem
(\ref{pbm:h-eq}), (\ref{pbm:h-data}) admits infinitely many local
solutions.  The second result is a detailed study of the case where
$u_{0}$ and $u_{1}$ are \emph{eigenvectors} of $A$ relative to the
\emph{same eigenvalue}.  In this very special case the authors proved
that uniqueness of the local solution fails if and only if the
following three conditions are satisfied:
\begin{list}{}{\leftmargin 4em \labelwidth 4em}
	\item[(AS1)] $\langle Au_{0},u_{1}\rangle= 0$;

	\item[(AS2)] $|A^{1/2}u_{1}|^{2}-\m{u_{0}}|Au_{0}|^{2}=0$;

	\item[(AS3)]  $m$ satisfies a suitable integrability condition in a
	neighborhood of $|A^{1/2}u_{0}|^{2}$.
\end{list}

In particular the local solution is unique if at least one of the
conditions above is not satisfied.

In this paper we extend the first two parts of this result to the
general case.  In Theorem~\ref{thm:main} we prove indeed that if
either condition (AS1) or condition (AS2) is not satisfied, then even
in the general case the local solution is always unique.

The proof of this result relies on two main steps.  

The first step is what we call \emph{trajectory uniqueness}.  We prove
indeed that the image of the curve $(u(t),u'(t))$ in the phase space
is unique.  To this end we parametrize the curve using the variable
$s=|A^{1/2}u(t)|^{2}$ instead of the variable $t$.  In this new
variable the trajectory is the image of a curve $(z(s),w(s))$, where
$z(s)$ and $w(s)$ are the solutions of a system in which the
non-Lipschitz nonlinear term $\m{u(t)}$ has become a non-Lipschitz
coefficient $m(s)$, which doesn't affect uniqueness.

The second step is what we call \emph{parametrization uniqueness}.
We prove indeed that the unique trajectory obtained in the first step
can be described by the solutions in a unique way.  To this end we
show that the parametrization $s(t)=|A^{1/2}u(t)|^{2}$ satisfies a
first order autonomous ordinary differential equation with
non-Lipschitz right-hand side, to which we can apply a uniqueness
result for nonstationary solutions (Lemma~\ref{lemma:ode-uniq}).

This paper is organized as follows.  In section~\ref{sec:statements}
we recall the definition of continuity modulus and Gevrey-type
functional spaces.  Moreover we state the classical local existence
result for (\ref{pbm:h-eq}), (\ref{pbm:h-data})
(Theorem~\ref{thm:hirosawa}) and our uniqueness result
(Theorem~\ref{thm:main}).  In section~\ref{sec:proofs} we prove
Theorem~\ref{thm:main}.  In section~\ref{sec:open} we collect some
open problems concerning uniqueness of solutions.

\setcounter{equation}{0}
\section{Preliminaries and statements}\label{sec:statements}

For the sake of simplicity we assume that $H$ admits a countable
complete orthonormal system $\{e_{k}\}_{k\geq 1}$ made by eigenvectors
of $A$.  We denote the corresponding eigenvalues by $\lambda_{k}^{2}$
(with $\lk\geq 0$), so that $Ae_{k}=\lambda_{k}^{2}e_{k}$ for every
$k\geq 1$.

Under this assumption we can work with Fourier series.  However, any
definition or statement of this section can be easily extended to the
general setting just by using the spectral decomposition instead of
Fourier series.  The interested reader is referred to \cite{ap} for
further details.

By means of the orthonormal system every $u\in H$ can be written in a
unique way in the form $u=\sum_{k=1}^{\infty}u_{k}e_{k}$, where
$u_{k}=\langle u,e_{k}\rangle$ are the Fourier components of $u$.
With these notations for every $\alpha\geq 0$ we have that 
$$D(A^{\alpha}):=\left\{u\in H:\sum_{k=1}^{\infty}
\lambda_{k}^{4\alpha}u_{k}^{2}<+\infty\right\}.$$

Let now $\varphi:[0,+\infty)\to(0,+\infty)$ be any function.  Then for
every $\alpha\geq 0$ and $r>0$ one can set
\begin{equation}
	\trebar{u}_{\varphi,r,\alpha}^{2}:=\sum_{k=1}^{\infty}\lambda_{k}^{4\alpha}
	u_{k}^{2} \exp\left(\strut r\varphi(\lambda_{k})\right),
	\label{defn:trebar}
\end{equation}
and then define the spaces
$$\G_{\varphi,r,\alpha}(A):=
\left\{u\in H:\trebar{u}_{\varphi,r,\alpha}<+\infty\right\}.$$

These spaces are a generalization of the usual spaces of Sobolev,
Gevrey or analytic functions.  They are Hilbert spaces with norm
$(|u|^{2}+\trebar{u}_{\varphi,r,\alpha}^{2})^{1/2}$.

A \emph{continuity modulus} is a continuous increasing function
$\omega:[0,+\infty)\to[0,+\infty)$ such that $\omega(0)=0$, and
$\omega(a+b)\leq\omega(a)+\omega(b)$ for every $a\geq 0$ and $b\geq
0$.

The function $m$ is said to be $\omega$-continuous if there exists a
constant $L\in\re$ such that
\begin{equation}
	|m(a)-m(b)|\leq
	L\,\omega(|a-b|)
	\hspace{3em}
	\forall a\geq 0,\ \forall b\geq 0.
	\label{hp:ocont}
\end{equation}

The following result sums up the state of the art concerning existence
of local solutions (see Theorem~2.1 and Theorem~2.2
in~\cite{hirosawa-main}, and the counterexamples in
\cite{gg:k-derloss}).

\begin{thmbibl}\label{thm:hirosawa}
	Let $\omega$ be a continuity modulus, let
	$m:[0,+\infty)\to[0,+\infty)$ be an \ocont\ function, and let
	$\varphi:[0,+\infty)\to(0,+\infty)$.
	
	Let us assume that there exists a constant $\Lambda$ such that
	\begin{equation}
		\sigma
		\omega\left(\frac{1}{\sigma}\right)\leq\Lambda\varphi(\sigma)
		\quad\quad
		\forall\sigma> 0
		\label{hp:phi-ndg}
	\end{equation}
	in the strictly hyperbolic case, and 
	\begin{equation}
		\sigma\leq\Lambda\varphi\left(\frac{\sigma}{
		\sqrt{\omega(1/\sigma)}}\right)
		\quad\quad
		\forall\sigma> 0
		\label{hp:phi-dg}
	\end{equation}
	in the weakly hyperbolic case.
	
	Let $(u_{0},u_{1})\in
	\G_{\varphi,r_{0},3/4}(A)\times\G_{\varphi,r_{0},1/4}(A)$ for some
	$r_{0}>0$.
	
	Then there exists $T>0$, and a nonincreasing function
	$r:[0,T]\to(0,r_{0}]$ such that problem (\ref{pbm:h-eq}),
	(\ref{pbm:h-data}) admits at least one local solution 
	\begin{equation}
		u\in C^{1}\left([0,T];\G_{\varphi,r(t),1/4}(A)\right)\cap
		C^{0}\left([0,T];\G_{\varphi,r(t),3/4}(A)\right).  
	\label{th:reg-sol}
\end{equation}
\end{thmbibl}

The main result of this paper is the following uniqueness result
for these solutions.

\begin{thm}\label{thm:main}
	Let $\omega$, $m$, $\varphi$ be as in Theorem~\ref{thm:hirosawa}.
	Let us assume that
	\begin{equation}
		(u_{0},u_{1})\in\G_{\varphi,r_{0},3/2}(A)\times
		\G_{\varphi,r_{0},1}(A) 
		\label{hp:main-data}
	\end{equation}
	for some $r_{0}>0$, and
	\begin{equation}
		\left|\langle Au_{0},u_{1}\rangle\right|+
		\left||A^{1/2}u_{1}|^{2}-\m{u_{0}}|Au_{0}|^{2}\right|\neq 0.
		\label{hp:main}
	\end{equation}
	
	Let us assume that problem (\ref{pbm:h-eq}), (\ref{pbm:h-data})
	admits two local solutions $v_{1}$ and $v_{2}$ in
	\begin{equation}
		C^{2}\left([0,T];\G_{\varphi,r(t),1/2}(A)\right)\cap
		C^{1}\left([0,T];\G_{\varphi,r(t),1}(A)\right)\cap
		C^{0}\left([0,T];\G_{\varphi,r(t),3/2}(A)\right)
		\label{hp:reg-sol}
	\end{equation}
	for some $T>0$, and some nonincreasing function
	$r:[0,T]\to(0,r_{0}]$.
	
	Then we have the following conclusions.
	\begin{enumerate}
		\renewcommand{\labelenumi}{(\arabic{enumi})}
		\item  There exists $T_{1}\in(0,T]$ such that
		\begin{equation}
			v_{1}(t)=v_{2}(t)
			\hspace{2em}
			\forall t\in[0,T_{1}].
			\label{th:uniqueness}
		\end{equation}
	
		\item  Let $T_{*}$ denote the supremum of all $T_{1}\in (0,T]$
		for which (\ref{th:uniqueness}) holds true. Let $v(t)$ denote 
		the common value of $v_{1}$ and $v_{2}$ in $[0,T_{*}]$.
		
		Then either $T_{*}=T$ or
		\begin{equation}
			\left|\langle Av(T_{*}),v'(T_{*})\rangle\right|+
			\left||A^{1/2}v'(T_{*})|^{2}-
			\m{v(T_{*})}|Av(T_{*})|^{2}\right|= 0.
			\label{th:limit}
		\end{equation}
	\end{enumerate}
\end{thm}

\begin{rmk}\label{rmk:reg-sol}
	\begin{em}
		The space (\ref{hp:reg-sol}) is the natural one when initial
		data satisfy (\ref{hp:main-data}).  Indeed from the linear
		theory it follows that any solution $u(t)$ of (\ref{pbm:h-eq})
		with
		$$u\in C^{0}([0,T];D(A^{3/4}))\cap C^{1}([0,T];D(A^{1/4}))$$
		and initial data as in (\ref{hp:main-data}) lies actually in
		(\ref{hp:reg-sol}).
	\end{em}
\end{rmk}
		
\begin{rmk}
	\begin{em}		
		Assumption (\ref{hp:main-data}) on the initial data is
		stronger than the corresponding assumption in
		Theorem~\ref{thm:hirosawa}.  This is due to a technical point
		in the proof.
		
		However in most cases the difference is only apparent. For
		example if $\omega(\sigma)=\sigma^{\beta}$ for some
		$\beta\in(0,1]$, then the following implication 
		$$u\in\G_{\varphi,r,0}(A)\ \Longrightarrow\ 
		u\in\G_{\varphi,r-\ep,\alpha}(A)$$
		holds true for every $r>0$, $\ep\in(0,r)$, $\alpha\geq 0$.
		Therefore in this case every solution satisfying
		(\ref{th:reg-sol}) fulfils (\ref{hp:reg-sol}) up to replacing 
		$r(t)$ with $r(t)/2$.
	\end{em}
\end{rmk}

\setcounter{equation}{0}
\section{Proofs}\label{sec:proofs}

\subsection{Technical lemmata}

\begin{lemma}\label{lemma:omega}
	Let $\omega:[0,+\infty)\to[0,+\infty)$ be a continuity modulus.
	
	Then 
	\begin{equation}
		\omega(x)\geq \omega(1)\cdot\frac{x}{x+1}
			\quad\quad
			\forall x\geq 0.
		\label{th:omega}
	\end{equation}
\end{lemma}

\prf Inequality (\ref{th:omega}) is trivial for $x=0$.  From the
subadditivity of $\omega$ it follows that $\omega(\lambda
x)\leq(\lambda+1)\omega(x)$ for every $\lambda\geq 0$ and $x\geq 0$
(this can be easily proved by induction on the integer part of
$\lambda$).  Applying this inequality with $x>0$ and
$\lambda=1/x$ we obtain (\ref{th:omega}) for $x>0$. 
\qed

\begin{lemma}\label{lemma:ode-lin}
	For $i=1,2$ let $\eta_{i}:(0,T]\to[0,+\infty)$ be a continuous
	function with finite integral. Let
	$y\in C^{0}([0,T];\re)\cap C^{1}((0,T];\re)$ be a function 
	such that $y(0)=0$, and
	\begin{equation}
		y'(t)\leq \eta_{1}(t)y(t)+\eta_{2}(t)
		\quad\quad
		\forall t\in(0,T].
		\label{hp:subsol}
	\end{equation}
	
	Then
	\begin{equation}
		y(t)\leq\exp\left(\int_{0}^{t}\eta_{1}(\tau)d\tau\right)\cdot
		\int_{0}^{t}\eta_{2}(\tau)\,d\tau
		\quad\quad
		\forall t\in[0,T].
		\label{th:supsol}
	\end{equation}
\end{lemma}

\prf
Let us consider the ordinary differential equation
\begin{equation}
	v'(t)=\eta_{1}(t)v(t)+\eta_{2}(t).
	\label{eq:ode}
\end{equation}

Assumption (\ref{hp:subsol}) is equivalent to say that $y$ is a
subsolution of (\ref{eq:ode}).  Since $\eta_{1}(t)$ and $\eta_{2}(t)$
are nonnegative it is easy to verify that the right-hand side of
(\ref{th:supsol}) is a supersolution of (\ref{eq:ode}).  Therefore
estimate (\ref{th:supsol}) follows from the standard comparison
principle.  \qed

\begin{lemma}\label{lemma:ode2ord}
	Let $y:[0,T]\to[0,+\infty)$ be a continuous function.  Let us
	assume that there exists $k\geq 0$ such that
	$$y(t)\leq k\int_{0}^{t}\frac{1}{s\sqrt{s}}
	\int_{0}^{s}y(\sigma)\,d\sigma\,ds.$$
	
	Then $y(t)=0$ for every $t\in[0,T]$.
\end{lemma}

\prf
Let us set $M:=\max\{y(t):t\in[0,T]\}$. Then an easy induction gives
$$y(t)\leq\frac{4^{n}k^{n}M}{n!}t^{n/2}
\quad\quad
\forall t\in[0,T],\quad
\forall n\in\n,$$
which implies the conclusion.
\qed

\begin{lemma}\label{lemma:ode-uniq}
	Let $s_{0}>0$, let $g:[0,s_{0}]\to\re$ be a continuous
	function, and let $T>0$.
	
	Then there exists at most one function $y:[0,T]\to[0,s_{0}]$ such 
	that
	\begin{eqnarray}
		 & y(0)=0, & 
		\label{hp:y-data} \\
		 & y'(t)>0
		\quad\quad
		\forall t\in(0,T], & 
		\label{hp:y'-pos}  \\
		 & y'(t)=g(y(t))
		\quad\quad
		\forall t\in(0,T]. & 
		\label{hp:y'=g}
	\end{eqnarray}
\end{lemma}

\prf Let $y_{1}(t)$ and $y_{2}(t)$ be two solutions of
(\ref{hp:y-data}), (\ref{hp:y'-pos}), (\ref{hp:y'=g}).  Let
$s_{1}:=y_{1}(T)$, $s_{2}:=y_{2}(T)$.  By (\ref{hp:y'-pos}) the
functions $y_{1}:[0,T]\to[0,s_{1}]$ and $y_{2}:[0,T]\to[0,s_{2}]$ are
strictly increasing and invertible.  Their inverse functions
$z_{1}(s)$ and $z_{2}(s)$ are defined and continuous in $[0,s_{3}]$,
where $s_{3}:=\min\{s_{1},s_{2}\}>0$.

Moreover $z_{1}$ and $z_{2}$ are of class $C^{1}$ in $(0,s_{3}]$, and
by (\ref{hp:y'=g}) 
$$z_{1}'(s)-z_{2}'(s)=
\frac{1}{y_{1}'(z_{1}(s))}-\frac{1}{y_{2}'(z_{2}(s))}=
\frac{1}{g(s)}-\frac{1}{g(s)}=0 \quad\quad
\forall s\in(0,s_{3}].$$

Since by (\ref{hp:y-data}) we have that $z_{1}(0)=z_{2}(0)=0$, it
follows that $z_{1}(s)=z_{2}(s)$ for every $s\in(0,s_{3}]$, and in
particular $s_{1}=s_{2}=y_{1}(T)=y_{2}(T)$.

Therefore also the inverse functions of $z_{1}$ and $z_{2}$, namely
$y_{1}$ and $y_{2}$, coincide.
\qed

\subsection{A variable change}

Let $u(t)$ be any solution of (\ref{pbm:h-eq}) defined in an interval
$[0,T]$.  Let us assume that $u$ belongs to the space
(\ref{hp:reg-sol}), and its initial data (\ref{pbm:h-data}) satisfy
(\ref{hp:main}).  Let us set
\begin{equation}
	\psi(t):=|A^{1/2}u(t)|^{2}-|A^{1/2}u_{0}|^{2}.
	\label{defn:psi}
\end{equation}

Then $\psi\in C^{2}([0,T])$, and 
$$\psi(0)=0,
\hspace{2em}
\psi'(0)=2\langle Au_{0},u_{1}\rangle,
\hspace{2em}
\psi''(0)=2\left(|A^{1/2}u_{1}|^{2}-\m{u_{0}}|Au_{0}|^{2}\right).$$

Our assumption (\ref{hp:main}) is equivalent to say that either
$\psi'(0)\neq 0$ or $\psi''(0)\neq 0$.  In both cases we can conclude
that there exists $T_{0}\in(0,T]$ such that $\psi'(t)$ has constant
sign in the interval $(0,T_{0}]$.

Let us assume, without loss of generality, that $\psi'(t)>0$ in
$(0,T_{0}]$.  Setting $s_{0}=\psi(T_{0})$, this implies that
$\psi:[0,T_{0}]\to[0,s_{0}]$ is strictly increasing and invertible.
Its inverse function $\psi^{-1}: [0,s_{0}]\to[0,T_{0}]$ belongs to
$C^{0}([0,s_{0}])\cap C^{2}((0,s_{0}])$, and
\begin{equation}
	(\psi^{-1})'(s)=\frac{1}{\psi'(\psi^{-1}(s))}= \frac{1}{2\langle
	Au(\psi^{-1}(s)),u'(\psi^{-1}(s))\rangle}>0 \quad\quad
	\forall s\in(0,s_{0}].
	\label{eq:phi'}
\end{equation}

Let us set now
\begin{equation}
	z(s):=A^{1/2}u(\psi^{-1}(s)),
	\hspace{3em}
	w(s):=u'(\psi^{-1}(s)).
	\label{defn:zw}
\end{equation}

From the regularity of $u$ and $\psi^{-1}$ it follows that $z(s)$ and 
$w(s)$ belong to
\begin{equation}
	C^{0}\left([0,s_{0}],\G_{\varphi,r_{1},1}\right)\cap
	C^{1}\left((0,s_{0}],\G_{\varphi,r_{1},1/2}\right)
	\label{defn:r1}
\end{equation}
for some $r_{1}>0$.  Moreover they satisfy the
initial conditions
\begin{equation}
	z(0)=A^{1/2}u_{0},
	\hspace{3em}
	w(0)=u_{1}.
	\label{eq:zw-data}
\end{equation}

The derivatives of $z(s)$ and $w(s)$ with respect to the variable $s$
can be easily computed using (\ref{pbm:h-eq}) and (\ref{eq:phi'}).
For every $s\in(0,s_{0}]$ it turns out that
\begin{eqnarray}
	z'(s) & = & \frac{A^{1/2}w(s)}{2\langle A^{1/2}z(s),w(s)\rangle},
	\label{eq:z}  \\
	w'(s) & = & -c(s)\frac{A^{1/2}z(s)}{2\langle
	A^{1/2}z(s),w(s)\rangle},
	\label{eq:w}
\end{eqnarray}
where $c(s):=m(s+|A^{1/2}u_{0}|^{2})$.

This system is singular when denominators vanish for $s=0$, i.e., when
$\langle Au_{0},u_{1}\rangle=0$.  However we claim that there exists
$s_{1}\in(0,s_{0}]$ such that ($\gamma_{1}$ is the first of a long
list of constants)
\begin{equation}
	\langle A^{1/2}z(s),w(s)\rangle\geq\gamma_{1}\sqrt{s}
	\quad\quad
	\forall s\in(0,s_{1}].
	\label{est:denom}
\end{equation}

To this end we first remark that
\begin{equation}
	\frac{\mathrm{d}}{\mathrm{d}s}\left(\langle
	A^{1/2}z,w\rangle^{2}\right)=
	|A^{1/2}w(s)|^{2}-c(s)|A^{1/2}z(s)|^{2},
	\label{deriv-denom}
\end{equation}
hence (we recall that $\psi'$ is assumed to be positive) 
\begin{equation}
	\langle A^{1/2}z(s),w(s)\rangle=\left[
	\langle Au_{0},u_{1}\rangle^{2}+
	\int_{0}^{s}\left(|A^{1/2}w(\sigma)|^{2}-
	c(\sigma)|A^{1/2}z(\sigma)|^{2}\right)
	\,d\sigma\right]^{1/2}.
	\label{eq:denom-int}
\end{equation}

If $\langle Au_{0},u_{1}\rangle>0$, then (\ref{est:denom}) is trivial
provided that $s_{1}$ is small enough.  If $\langle
Au_{0},u_{1}\rangle=0$, then assumption (\ref{hp:main}) implies that
$|A^{1/2}u_{1}|^{2}-\m{u_{0}}|Au_{0}|^{2}>0$, hence the right-hand
side of (\ref{deriv-denom}) is larger than a positive constant in a
right neighborhood of $0$, so that (\ref{est:denom}) follows from
(\ref{eq:denom-int}).

\subsection{Trajectory uniqueness}

Let $v_{1}(t)$ and $v_{2}(t)$ be two solutions of (\ref{pbm:h-eq}),
(\ref{pbm:h-data}).  Let us define $\psi_{1}(t)$ and $\psi_{2}(t)$
according to (\ref{defn:psi}), and then $(z_{1}(s),w_{1}(s))$ and
$(z_{2}(s),w_{2}(s))$ according to (\ref{defn:zw}).  Let $s_{1}>0$ be
small enough so that $z_{1}(s)$, $z_{2}(s)$, $w_{1}(s)$, $w_{2}(s)$
are defined in $[0,s_{1}]$, and in this interval they are as regular
as prescribed by (\ref{defn:r1}), and they satisfy system
(\ref{eq:z}), (\ref{eq:w}), and estimate (\ref{est:denom}).

We claim that $z_{1}(s)=z_{2}(s)$ and $w_{1}(s)=w_{2}(s)$ in
$[0,s_{2}]$ for a suitable $s_{2}\in(0,s_{1}]$.  To this end we
introduce the differences
\begin{equation}
	x(s):=z_{1}(s)-z_{2}(s),
	\hspace{3em}
	y(s):=w_{1}(s)-w_{2}(s).
	\label{defn:xy}
\end{equation}

Setting for simplicity
$$d_{1}(s):=2\langle A^{1/2}z_{1}(s),w_{1}(s)\rangle,
\hspace{3em}
d_{2}(s):=2\langle A^{1/2}z_{2}(s),w_{2}(s)\rangle,$$
it is easy to see that $x(s)$ and $y(s)$ are solutions in $(0,s_{1}]$
of the system
\begin{eqnarray}
	x'(s) & = & \frac{A^{1/2}y(s)}{d_{1}(s)} +\left(
	\frac{1}{d_{1}(s)}-\frac{1}{d_{2}(s)}\right)A^{1/2}w_{2}(s),
	\label{eq:x}  \\
	y'(s) & = & -c(s)\frac{A^{1/2}x(s)}{d_{1}(s)} -c(s)\left(
	\frac{1}{d_{1}(s)}-\frac{1}{d_{2}(s)}\right)A^{1/2}z_{2}(s),
	\label{eq:y}
\end{eqnarray}
with initial data $x(0)=y(0)=0$.

Let us introduce the Fourier components $x_{k}(s)$, $y_{k}(s)$,
$z_{i,k}(s)$, $w_{i,k}(s)$ of $x(s)$, $y(s)$, $z_{i}(s)$, $w_{i}(s)$
(with $i=1,2$).  System (\ref{eq:x}), (\ref{eq:y}) becomes a system of
infinitely many ordinary differential equations of the form
\begin{eqnarray}
	x_{k}'(s) & = & \frac{\lk y_{k}(s)}{d_{1}(s)} +\lk\left(
	\frac{1}{d_{1}(s)}-\frac{1}{d_{2}(s)}\right)w_{2,k}(s),
	\label{eq:xk}  \\
	y_{k}'(s) & = & -c(s)\frac{\lk x_{k}(s)}{d_{1}(s)} -c(s)\lk\left(
	\frac{1}{d_{1}(s)}-\frac{1}{d_{2}(s)}\right)z_{2,k}(s),
	\label{eq:yk}
\end{eqnarray}
all with initial data $x_{k}(0)=y_{k}(0)=0$.

If $\lk=0$ it is clear that $x_{k}(s)=y_{k}(s)=0$ in $[0,s_{1}]$.  So
let us concentrate on the components corresponding to positive
eigenvalues.  To this end we consider the \emph{approximated energy
estimates} introduced in~\cite{dgcs} and \cite{cjs}, which are
different in the strictly hyperbolic and in the weakly hyperbolic
case.

\paragraph{The strictly hyperbolic case}

Let us assume that
\begin{equation}
	m(\sigma)\geq\gamma_{2}>0
	\quad\quad
	\forall\sigma\geq 0.
	\label{est:ndg}
\end{equation}

In particular the same estimate holds true for $c(s)$.  Formally we
need $c(s)$ to be defined only for $s\in[0,s_{1}]$.  In order to make
convolutions we extend $c(s)$ to the whole real line by setting
$c(s)=c(0)$ for every $s\leq 0$, and $c(s)=c(s_{1})$ for every $s\geq
s_{1}$.

Let us fix once for all a function $\rho:\re\to [0,+\infty)$ of class
$C^{\infty}$, with compact support and integral equal to 1. For every 
$\ep>0$ let us set
$$\cep(s):=\int_{\re}^{}c(s+\ep\sigma)\rho(\sigma)\,d\sigma.$$

From the boundedness and the $\omega$-continuity of $c(s)$ it is easy
to deduce that for every $s\in[0,s_{1}]$ (actually for every
$s\in\re$) we have that (from now on all constants are independent on
$\ep$)
\begin{eqnarray}
	 & |\cep(s)-c(s)|\leq\gamma_{3}\omega(\ep), & 
	\label{est:cep-c}  \\
	 & |\cep'(s)|\leq\gamma_{4}\displaystyle{\frac{\omega(\ep)}{\ep}}, & 
	\label{est:cep'}  \\
	 & \gamma_{2}\leq\cep(s)\leq\gamma_{5}. & 
	\label{est:cep-nu}
\end{eqnarray}

Let us consider the energy
\begin{equation}
	E_{k,\ep}(s):=|y_{k}|^{2}+\cep(s)|x_{k}(s)|^{2}.
	\label{defn:ekep}
\end{equation}

From (\ref{eq:xk}) and (\ref{eq:yk}) we have that
\begin{eqnarray}
	E_{k,\ep}'(s) & = & \cep'(s)|x_{k}|^{2}+2(\cep(s)-c(s)) \frac{\lk
	x_{k}y_{k}}{d_{1}(s)}+ 
	\nonumber  \\
	 & & +2\lk\left(
	\frac{1}{d_{1}(s)}-\frac{1}{d_{2}(s)}\right)\left( \cep(s)
	x_{k}w_{2,k}-c(s) y_{k}z_{2,k}\right)
	\nonumber  \\
	 & =: & I_{1}(s)+I_{2}(s)+I_{3}(s).
	\label{eq:ekep'}
\end{eqnarray} 

Let us estimate the three terms.  By (\ref{est:cep'}) and
(\ref{est:cep-nu}) we have that
\begin{equation}
	I_{1}(s)\leq\frac{|\cep'(s)|}{\cep(s)}\cdot\cep(s)|x_{k}(s)|^{2}\leq
	\gamma_{6}\frac{\omega(\ep)}{\ep}E_{\ep,k}(s)\leq
	\gamma_{7}\frac{\omega(\ep)}{\ep}\cdot
	\frac{E_{\ep,k}(s)}{\sqrt{s}}.
	\label{est:i1}
\end{equation}

By (\ref{est:denom}), (\ref{est:cep-c}), and (\ref{est:cep-nu}) we
have that
\begin{equation}
	I_{2}(s)\leq 2\lk\frac{|\cep(s)-c(s)|}{d_{1}(s)\sqrt{\cep(s)}}\cdot
	|y_{k}(s)|\cdot\sqrt{\cep(s)}|x_{k}(s)|\leq
	\gamma_{8}\lk\frac{\omega(\ep)}{\sqrt{s}}E_{k,\ep}(s).
	\label{est:i2}
\end{equation}

It remains to estimate $I_{3}(s)$. Since the norms $|A^{1/2}z_{i}(s)|$
and $|A^{1/2}w_{i}(s)|$ are bounded we have that
$$\left||A^{1/2}z_{1}|^{2}-|A^{1/2}z_{2}|^{2}\right|=
\left|\langle A^{1/2}(z_{1}+z_{2}),A^{1/2}(z_{1}-z_{2})\rangle\right|
\leq\gamma_{9}|A^{1/2}x|,$$
$$\left||A^{1/2}w_{1}|^{2}-|A^{1/2}w_{2}|^{2}\right|=
\left|\langle A^{1/2}(w_{1}+w_{2}),A^{1/2}(w_{1}-w_{2})\rangle\right|
\leq\gamma_{10}|A^{1/2}y|,$$
hence by (\ref{deriv-denom}) and the boundedness of $c(s)$
$$\left|\frac{\mathrm{d}}{\mathrm{d}s}
\left(d_{1}^{2}(s)-d_{2}^{2}(s)\right)\right|\leq\gamma_{11}\left(
|A^{1/2}x(s)|+|A^{1/2}y(s)|\right).$$

It follows that
\begin{equation}
	\left|d_{1}^{2}(s)-d_{2}^{2}(s)\right|\leq
	\gamma_{11}\int_{0}^{s}\left(
	|A^{1/2}x(\sigma)|+|A^{1/2}y(\sigma)|\right)d\sigma=:
	\psi_{1,2}(s),
	\label{defn:psi12}
\end{equation}
hence by (\ref{est:denom})
$$\left|\frac{1}{d_{1}(s)}-\frac{1}{d_{2}(s)}\right|=
\frac{\left|d_{2}^{2}(s)-d_{1}^{2}(s)\right|}{d_{1}(s)d_{2}(s)
(d_{1}(s)+d_{2}(s))}\leq\gamma_{12}
\frac{1}{s\sqrt{s}}\psi_{1,2}(s).$$

Since $c(s)$ and $\cep(s)$ are bounded from above we have that
$$|\cep(s) x_{k}w_{2,k}-c(s) y_{k}z_{2,k}|\leq\gamma_{13}
\left(\sqrt{\cep(s)}|x_{k}|\cdot|w_{2,k}|+
|y_{k}|\cdot|z_{2,k}|\right),$$
hence
\begin{eqnarray}
	I_{3}(s) & \leq & \frac{\gamma_{14}}{\sqrt{s}}\left(
	\frac{\psi_{1,2}(s)}{s}\lk|w_{2,k}|\cdot
	|\sqrt{\cep(s)}x_{k}|+\frac{\psi_{1,2}(s)}{s}\lk|z_{2,k}|\cdot
	|y_{k}|\right)
	\nonumber  \\
	 & \leq & \frac{\gamma_{15}}{\sqrt{s}}\left(
	 \frac{\psi_{1,2}^{2}(s)}{s^{2}}\lk^{2}|w_{2,k}|^{2}+\cep(s)|x_{k}|^{2}+
	 \frac{\psi_{1,2}^{2}(s)}{s^{2}}\lk^{2}|z_{2,k}|^{2}+
	 |y_{k}|^{2}\right)
	\nonumber  \\
	 & = & \frac{\gamma_{15}}{\sqrt{s}}E_{k,\ep}+
	 \gamma_{15}\frac{\psi_{1,2}^{2}(s)}{s^{2}\sqrt{s}}\lk^{2}
	 \left(|w_{2,k}|^{2}+|z_{2,k}|^{2}\right).
	\label{est:i3}
\end{eqnarray}

From (\ref{eq:ekep'}), (\ref{est:i1}), (\ref{est:i2}), (\ref{est:i3}) 
we therefore obtain that
$$E_{k,\ep}'\leq\gamma_{16}\left( \frac{\omega(\ep)}{\ep}+
\lk\omega(\ep)+1\right) \frac{E_{k,\ep}}{\sqrt{s}}+\gamma_{15}
\frac{\psi_{1,2}^{2}(s)}{s^{2}\sqrt{s}}\lk^{2}\left(
|w_{2,k}|^{2}+|z_{2,k}|^{2}\right).$$

Let us set now $\ep_{k}=\lk^{-1}$ (we recall that we can limit
ourselves to positive eigenvalues).  By assumption (\ref{hp:phi-ndg})
we have that 
$$\frac{\omega(\ep_{k})}{\ep_{k}}=\lk\omega(\ep_{k})=
\lk\omega\left(\frac{1}{\lk}\right)\leq
\Lambda\varphi(\lk),$$
hence
\begin{eqnarray*}
	E_{k,\ep_{k}}'(s) & \leq &
	\gamma_{17}\frac{\varphi(\lk)+1}{\sqrt{s}}\cdot
	E_{k,\ep_{k}}(s)+\gamma_{15}
	\frac{\psi_{1,2}^{2}(s)}{s^{2}\sqrt{s}}\lk^{2}\left(
	|w_{2,k}(s)|^{2}+|z_{2,k}(s)|^{2}\right) \\
	 & =: & \eta_{1}(s)E_{k,\ep_{k}}(s)+\eta_{2}(s).
\end{eqnarray*}

The integral of $\eta_{1}(s)$ in $[0,s_{1}]$ is finite. Moreover from 
definition (\ref{defn:psi12}) of $\psi_{1,2}$ it is clear that
$\psi_{1,2}(s)\leq\gamma_{18}s$. It follows that also the 
integral of $\eta_{2}(s)$ in $[0,s_{2}]$ is finite. 

We can therefore apply Lemma~\ref{lemma:ode-lin}.  Since for every
$s\in[0,s_{1}]$ we have that
$$\exp\left(\int_{0}^{s}\eta_{1}(\sigma)\,d\sigma\right)=
\exp\left(2\gamma_{17}\varphi(\lk)\sqrt{s}+
2\gamma_{17}\sqrt{s}\right)\leq
\gamma_{19}\exp\left(\gamma_{20}\varphi(\lk)\sqrt{s}\right),$$
it follows that
\begin{equation}
	E_{k,\ep_{k}}(s)\leq\gamma_{21}\exp\left(
	\gamma_{20}\varphi(\lk)\sqrt{s}\right)\int_{0}^{s}
	\frac{\psi_{1,2}^{2}(\sigma)}{\sigma^{2}\sqrt{\sigma}}\lk^{2}\left(
	|z_{2,k}(\sigma)|^{2}+|w_{2,k}(\sigma)|^{2}\right)d\sigma.
	\label{est:ekek}
\end{equation}

Let us choose $s_{2}\in(0,s_{1}]$ such that
$\gamma_{20}\sqrt{s_{2}}\leq r_{1}$.  By (\ref{est:cep-nu}) and
(\ref{est:ekek}) we have that
\begin{eqnarray*}
	|y_{k}(s)|^{2}+|x_{k}(s)|^{2} & \leq & 
	\max\left\{1,\frac{1}{\gamma_{2}}\right\}
	E_{k,\ep_{k}}(s)\\
	 & \leq & \gamma_{22}\int_{0}^{s}
	 \frac{\psi_{1,2}^{2}(\sigma)}{\sigma^{2}\sqrt{\sigma}}
	 \lk^{2}\exp\left(r_{1}\varphi(\lk)\right)\left(
	|z_{2,k}(\sigma)|^{2}+|w_{2,k}(\sigma)|^{2}\right)
	 d\sigma.
\end{eqnarray*}

Summing over $k$ and recalling that $z_{2}$ and $w_{2}$ belong to the
space (\ref{defn:r1}) we find that
$$|A^{1/2}x(s)|^{2}+|A^{1/2}y(s)|^{2}=
\sum_{k=1}^{\infty}\lk^{2}\left(|x_{k}(s)|^{2}+|y_{k}(s)|^{2}\right)\leq$$
$$\leq\gamma_{22} \int_{0}^{s}
\frac{\psi_{1,2}^{2}(\sigma)}{\sigma^{2}\sqrt{\sigma}}\left(
\|z_{2}(\sigma)\|_{\varphi,r_{1},1}^{2}+
\|w_{2}(\sigma)\|_{\varphi,r_{1},1}^{2}\right)d\sigma
\leq\gamma_{23} \int_{0}^{s}
\frac{\psi_{1,2}^{2}(\sigma)}{\sigma^{2}\sqrt{\sigma}}\,d\sigma.$$

By definition (\ref{defn:psi12}) of $\psi_{1,2}$ and
H\"{o}lder's inequality we obtain that
\begin{eqnarray*}
	|A^{1/2}x(s)|^{2}+|A^{1/2}y(s)|^{2} & \leq & 
	\gamma_{23}\int_{0}^{s}\frac{1}{\sigma^{2}\sqrt{\sigma}}
	\left[\int_{0}^{\sigma}\left(|A^{1/2}x(\tau)|+
	|A^{1/2}y(\tau)|\right)d\tau\right]^{2}d\sigma\\
	 & \leq & \gamma_{23}\int_{0}^{s}\frac{1}{\sigma\sqrt{\sigma}}
	\int_{0}^{\sigma}\left(|A^{1/2}x(\tau)|^{2}+
	|A^{1/2}y(\tau)|^{2}\right)d\tau d\sigma.
\end{eqnarray*}

Applying Lemma~\ref{lemma:ode2ord} we conclude that
$|A^{1/2}x(s)|^{2}= |A^{1/2}y(s)|^{2}=0$ for every $s\in[0,s_{2}]$,
namely $z_{1}(s)=z_{2}(s)$ and $w_{1}(s)=w_{2}(s)$ in the same
interval.

\paragraph{The weakly hyperbolic case}

Let us modify $c(s)$ outside the interval $[0,s_{1}]$ as in the
strictly hyperbolic case. Since $c(s)$ is bounded we can also assume
that $\omega$ is bounded. For every $\ep>0$ let us set
$$\cep(s):=\omega(\ep)+\int_{\re}^{}c(s+\ep\sigma)\rho(\sigma)\,d\sigma.$$

Estimates (\ref{est:cep-c}) and (\ref{est:cep'}) are still true, but
(\ref{est:cep-nu}) has to be replaced by the weaker (for the estimate 
from above we need the boundedness of $\omega$)
\begin{equation}
	\omega(\ep)\leq\cep(s)\leq\gamma_{24}.
	\label{est:cep-nu-w}
\end{equation} 

Let us define $E_{k,\ep}(s)$ according to (\ref{defn:ekep}).  Its
derivative is the same as in the strictly hyperbolic case.  So we need
to estimates the three summands in (\ref{eq:ekep'}).  Using
(\ref{est:cep-nu-w}) instead of (\ref{est:cep-nu}) we find that
$$I_{1}(s)\leq\gamma_{25}\frac{1}{\ep}\cdot\frac{E_{k,\ep}(s)}{\sqrt{s}},
\hspace{3em}
I_{2}(s)\leq\gamma_{26}\lk\sqrt{\omega(\ep)}\cdot
\frac{E_{k,\ep}(s)}{\sqrt{s}}.$$

The estimate on $I_{3}(s)$ is exactly the same as in (\ref{est:i3}).
We finally obtain that
$$E_{k,\ep}'(s)\leq\gamma_{27}\left( \frac{1}{\ep}+
\lk\sqrt{\omega(\ep)}+1\right)
\frac{E_{k,\ep}(s)}{\sqrt{s}}+\gamma_{28}
\frac{\psi_{1,2}^{2}(s)}{s^{2}\sqrt{s}}\lk^{2}\left(
|w_{2,k}(s)|^{2}+|z_{2,k}(s)|^{2}\right).$$

Now we choose $\ep$ as a function of $k$.  We consider the function
$h(\ep):=\ep\sqrt{\omega(\ep)}$, which is invertible, and we set
$\ep_{k}:=h^{-1}(1/\lk)$.

Applying assumption (\ref{hp:phi-dg}) with $\sigma=1/\ep_{k}$ we
obtain that
\begin{equation}
	\frac{1}{\ep_{k}}\leq\Lambda\varphi\left(
	\frac{1}{\ep_{k}\sqrt{\omega(\ep_{k})}}\right)
	=\Lambda\varphi\left(\frac{1}{h(\ep_{k})}\right)
	=\Lambda\varphi\left(\lk\right),
	\label{est:eklk}
\end{equation}
and therefore
$$\frac{1}{\ep_{k}}+\lk\sqrt{\omega(\ep_{k})}=
\frac{1+h(\ep_{k})\lk}{\ep_{k}}= \frac{2}{\ep_{k}}\leq 2\Lambda
\varphi\left(\lk\right)=\gamma_{29}
\varphi(\lk),$$
hence
$$E_{k,\ep_{k}}'(s)\leq\gamma_{30}\frac{\varphi(\lk)+1}{\sqrt{s}}
E_{k,\ep_{k}}(s)+\gamma_{28}
\frac{\psi_{1,2}^{2}(s)}{s^{2}\sqrt{s}}\lk^{2}\left(
|w_{2,k}(s)|^{2}+|z_{2,k}(s)|^{2}\right).$$

As in the strictly hyperbolic case we can apply Lemma~\ref{lemma:ode-lin} 
to this differential inequality and obtain that
$$E_{k,\ep_{k}}(s)\leq\gamma_{31}\exp\left(
\gamma_{32}\varphi(\lk)\sqrt{s}\right)\int_{0}^{s}
\frac{\psi_{1,2}^{2}(\sigma)}{\sigma^{2}\sqrt{\sigma}}\lk^{2}
\left( |z_{2,k}(\sigma)|^{2}+|w_{2,k}(\sigma)|^{2}\right)d\sigma.$$

Let us choose $s_{2}\in(0,s_{1}]$ such that
$\gamma_{32}\sqrt{s_{2}}\leq r_{1}/2$.

Applying Lemma~\ref{lemma:omega} and (\ref{est:eklk}) we have that
$$\max\left\{1,\frac{1}{\omega(\ep_{k})}\right\}\leq
1+\frac{1}{\omega(\ep_{k})}\leq
\gamma_{33}\left(1+\frac{1}{\ep_{k}}\right)\leq
\gamma_{34}(\varphi(\lk)+1)\leq
\gamma_{35}\exp(r_{1}\varphi(\lk)/2),$$
independently on $k$, hence
\begin{eqnarray*}
	|y_{k}(s)|^{2}+|x_{k}(s)|^{2} & \leq & 
	\max\left\{1,\frac{1}{\omega(\ep_{k})}\right\}
	E_{k,\ep_{k}}(s)\\
	 & \leq & \gamma_{36}\int_{0}^{s}
	 \frac{\psi_{1,2}^{2}(\sigma)}{\sigma^{2}\sqrt{\sigma}}
	 \lk^{2}\exp\left(r_{1}\varphi(\lk)\right)\left(
	|z_{2,k}(\sigma)|^{2}+|w_{2,k}(\sigma)|^{2}\right)
	 d\sigma.
\end{eqnarray*}

From now on we proceed exactly as in the strictly hyperbolic case.

\subsection{Parametrization uniqueness}

Let us come back to the two solutions $v_{1}(t)$ and $v_{2}(t)$ of
problem (\ref{pbm:h-eq}), (\ref{pbm:h-data}).  We already defined
$\psi_{1}(t)$ and $\psi_{2}(t)$ according to (\ref{defn:psi}), and
then $(z_{1},w_{1})$, and $(z_{2},w_{2})$ according to (\ref{defn:zw}).
For $i=1,2$ we have that
$$\psi_{i}'(t)=2\langle Av_{i}(t),v_{i}'(t)\rangle= 2\langle
Av_{i}(\psi_{i}^{-1}(\psi_{i}(t))),
v_{i}'(\psi_{i}^{-1}(\psi_{i}(t)))\rangle=$$
$$=2\langle A^{1/2}z_{i}(\psi_{i}(t)),
w_{i}(\psi_{i}(t))\rangle$$
for every small enough $t$.  Since $z_{1}(s)=z_{2}(s)=:z(s)$ and
$w_{1}(s)=w_{2}(s)=:w(s)$ in an interval $[0,s_{2}]$, we have that in
an interval $[0,T_{1}]$ the functions $\psi_{1}(t)$ and $\psi_{2}(t)$
are solutions of the Cauchy problem
$$\psi'(t)=2\langle A^{1/2}z(\psi(t)), w(\psi(t))\rangle=:g(\psi(t)),
\hspace{3em} 
\psi(0)=\langle Au_{0},u_{1}\rangle.$$

Since we already know that these solutions are strictly increasing in
$[0,T_{1}]$ we can apply Lemma~\ref{lemma:ode-uniq} and deduce that
$\psi_{1}(t)=\psi_{2}(t)$ in $[0,T_{1}]$. Finally we have that
$$v_{1}'(t)=v_{1}'(\psi_{1}^{-1}(\psi_{1}(t)))=
w_{1}(\psi_{1}(t))=w_{2}(\psi_{2}(t))=
v_{2}'(\psi_{2}^{-1}(\psi_{2}(t)))=v_{2}'(t)$$
in $[0,T_{1}]$, hence also $v_{1}(t)=v_{2}(t)$ in the same interval.

\subsection{Continuation}

Let us prove the second statement of Theorem~\ref{thm:main}. The
argument is quite standard. Let us
assume by contradiction that two solutions $v_{1}(t)$ and
$v_{2}(t)$ are defined in an interval $[0,T]$, and coincide in a
maximal interval $[0,T_{*}]$ with $T_{*}<T$.  If (\ref{th:limit}) is
not satisfied, then we can apply the first statement with ``initial''
data in $T_{*}$, and deduce that $v_{1}$ and $v_{2}$ coincide in some 
interval $[T_{*},T_{*}+\delta]$.

This contradicts the maximality of $T_{*}$.
\qed

\setcounter{equation}{0}
\section{Open problems}\label{sec:open}

The uniqueness problem for Kirchhoff equations is quite open. In this 
section we state four questions in this field.

The first one concerns counterexamples. We don't know any example
where uniqueness fails apart from those given in~\cite{as}. So we 
ask whether different counterexamples can be provided.

\begin{open}
	Let $\omega$, $m$, $\varphi$, $u_{0}$, $u_{1}$ be as in
	Theorem~\ref{thm:main}, but without assumption (\ref{hp:main}).
	Let us assume that problem (\ref{pbm:h-eq}), (\ref{pbm:h-data})
	admits two local solutions.
	
	Can we conclude that $u_{0}$ and $u_{1}$ are eigenvectors of $A$
	relative to the same eigenvalue?
\end{open}

We point out that this problem is open even in the simple case
$H=\re^{2}$, where $\omega$ and $\varphi$ play non role, and no
regularity is required on initial data.

The second open problem concerns trajectory uniqueness (the key step in
our proof).

\begin{open}
	Let $\omega$, $m$, $\varphi$, $u_{0}$, $u_{1}$ be as in
	Theorem~\ref{thm:main}, but without assumption (\ref{hp:main}).
	Let us consider system (\ref{eq:z}), (\ref{eq:w}), with initial
	data (\ref{eq:zw-data}).
	
	Does this system admit at most one solution?
\end{open}

Note that in the case where $\langle Au_{0},u_{1}\rangle=0$ it is by
no means clear that a solution always exists, since this implicitly
requires that $\langle A^{1/2}z(s),w(s)\rangle\neq 0$ for every
$s\in(0,s_{0}]$.  We point out that, even in the nonuniqueness
examples of \cite{as}, the solution of this system exists and it is
unique.

The third open problem concerns the regularity of \emph{initial data}.
It may happen indeed that problem (\ref{pbm:h-eq}), (\ref{pbm:h-data})
has a solution even for some initial data that do not satisfy the
assumptions of Theorem~\ref{thm:hirosawa} (see for example the
solutions with derivative loss constructed in \cite{gg:k-derloss}).
Are there uniqueness results for these solutions?

\begin{open}
	Is it possible to prove the known uniqueness results (namely the
	Lipschitz case and our Theorem~\ref{thm:main}) with less
	regularity requirements on initial data?
\end{open}

The last open problem concerns regularity of \emph{solutions}.  Both
the result in the Lipschitz case, and our result require the a priori
assumption that solutions lie in $D(A^{3/4})\times D(A^{1/4})$ (see
Remark~\ref{rmk:reg-sol}).  By the linear theory these solutions
automatically belong to the same space (technically to the same scale
of spaces) of the initial data.  On the other hand, equation
(\ref{pbm:h-eq}) makes perfectly sense in the energy space
$D(A^{1/2})\times H$.  Just to give an extreme example, let us
consider the strictly hyperbolic case, with a Lipschitz continuous
nonlinearity $m$, and analytic initial data.  We know that there is a
unique solution in $D(A^{3/4})\times D(A^{1/4})$, which is actually
analytic.  However as far as we know no one can exclude that there
exists a different solution in $D(A^{1/2})\times H$ with the same
(analytic) initial data.

\begin{open}
	Is it possible to extend the known uniqueness results to solutions
	in the energy space?
\end{open}

\label{NumeroPagine}

\end{document}